\documentclass[a4paper,12pt]{amsart}
\usepackage{amssymb}
\usepackage{ifthen}
\usepackage[dvips]{graphicx}
\nonstopmode \numberwithin{equation}{section}
\usepackage{amssymb}
\usepackage{ifthen}
\usepackage{graphicx}
\usepackage{amsmath}
\usepackage[T1]{fontenc} 
\usepackage[utf8]{inputenc}
\usepackage[usenames,dvipsnames]{color}
\usepackage{color}
\usepackage[english]{babel}
\usepackage{fancyhdr}
\usepackage{fancybox}
\usepackage{tikz}

\nonstopmode \numberwithin{equation}{section}
\theoremstyle{plain}

\newtheorem{conj}{Conjecture}

\theoremstyle{definition}
\newtheorem{defn}{Definition}[section]

\newtheorem{thm}{Theorem}[section]
\newtheorem{prob}{Problem}[section]
\newtheorem{cor}{Corollary}[section]
\newtheorem{ques}{Question}[section]
\newtheorem{prop}{Proposition}[section]
\newtheorem{rem}{Remark}[section]
\newtheorem{lem}{Lemma}[section]

\newcounter{minutes}
\divide\time by 60
\newcounter{hours}
\multiply\time by 60
\addtocounter{minutes}{-\time}

\newcounter {own}
\def\theown {\thesection       .\arabic{own}}

\newenvironment{pf}[1][]{%
	\vskip 3mm
	\noindent
	\ifthenelse{\equal{#1}{}}%
	{{\slshape Proof. }}%
	{{\slshape #1.} }%
}%
{\qed\bigskip}

\newcounter{alphabet}

\newcommand{\M}{{\mathcal M}}

\def\be{\begin{equation}}
	\def\ee{\end{equation}}

\newcommand{\bee}{\begin{enumerate}}
	\newcommand{\eee}{\end{enumerate}}

\newcommand{\blem}{\begin{lem}}
	\newcommand{\elem}{\end{lem}}
\newcommand{\bthm}{\begin{thm}}
	\newcommand{\ethm}{\end{thm}}
\newcommand{\bcor}{\begin{cor}}
	\newcommand{\ecor}{\end{cor}}
\newcommand{\beg}{\begin{examp}}
	\newcommand{\eeg}{\end{examp}}
\newcommand{\begs}{\begin{examples}}
	\newcommand{\eegs}{\end{examples}}

\newcommand{\bdefn}{\begin{defn}}
	\newcommand{\edefn}{\end{defn}}

\newcommand{\bprob}{\begin{prob}}
	\newcommand{\eprob}{\end{prob}}
\newcommand{\bei}{\begin{itemize}}
	\newcommand{\eei}{\end{itemize}}

\newcommand{\bcon}{\begin{conj}}
	\newcommand{\econ}{\end{conj}}
\newcommand{\bcons}{\begin{conjs}}
	\newcommand{\econs}{\end{conjs}}
\newcommand{\bprop}{\begin{prop}}
	\newcommand{\eprop}{\end{prop}}
\newcommand{\br}{\begin{rem}}
	\newcommand{\er}{\end{rem}}
\newcommand{\brs}{\begin{rems}}
	\newcommand{\ers}{\end{rems}}
\newcommand{\bo}{\begin{obser}}
	\newcommand{\eo}{\end{obser}}
\newcommand{\bos}{\begin{obsers}}
	\newcommand{\eos}{\end{obsers}}
\newcommand{\bpf}{\begin{pf}}
	\newcommand{\epf}{\end{pf}}
\newcommand{\ba}{\begin{array}}
	\newcommand{\ea}{\end{array}}
\newcommand{\beq}{\begin{eqnarray}}
	\newcommand{\beqq}{\begin{eqnarray*}}
		\newcommand{\eeq}{\end{eqnarray}}
	\newcommand{\eeqq}{\end{eqnarray*}}

\begin{document}
	
	\title{THE THIRD HANKEL DETERMINANT FOR INVERSE COEFFICIENTS OF STARLIKE FUNCTION OF ORDER 1/2}
	
	\author{Molla Basir Ahamed}
	\address{Molla Basir Ahamed, Department of Mathematics, Jadavpur University, Kolkata-700032, West Bengal, India.}
	\email{mbahamed.math@jadavpuruniversity.in}
	\author{Partha Pratim Roy}
	\address{Partha Pratim Roy, Department of Mathematics, Jadavpur University, Kolkata-700032, West Bengal, India.}
	\email{pproy.math.rs@jadavpuruniversity.in}
	
	\subjclass[{AMS} Subject Classification:]{Primary 30A10, 30H05, 30C35, Secondary 30C45}
	\keywords{Analytic functions, Starlike functions of order 1/2, Hankel determinant, inverse coefficients, invariance}
	
	\def\thefootnote{}
	\footnotetext{ {\tiny File:~\jobname.tex,
			printed: \number\year-\number\month-\number\day,
			\thehours.\ifnum\theminutes<10{0}\fi\theminutes }
	} \makeatletter\def\thefootnote{\@arabic\c@footnote}\makeatother
	
	\begin{abstract} 
		The sharp bound for the third Hankel determinant for the coefficients of the inverse function of starlike function of order  $1/2$ is obtained. In light of this, we can deduce that the functionals $|H_3(1)(f)|$ and $|H_3(1)(f^{-1})|$ exhibit invariance on the class $\mathcal{S}^*(1/2)$.
	\end{abstract}
	
	\maketitle
	\pagestyle{myheadings}
	\markboth{M. B. Ahamed and P. P. Roy}{The Third Hankel Determinant For Inverse Coefficients Of Starlike Function Of Order 1/2}
	
	\section{Introduction}
	In recent years, finding sharp bounds for Hankel determinants for various classes of univalent functions becomes a subject of great interest in geometric function theory. Establishing invariance between two functionals of Hankel determinants for a class of functions is an important study in Geometric function theory. There are results for certain classes of functions for Hankel determinant which are invariance (see e.g., \cite{Thomas-Verma-BAMS-2017}) and not invariance as well (see e.g., \cite{Raza-Riaz-Tho-BAMS-2023}). Thus it is understood that not every class of functions possesses this property concerning Hankel determinants. Therefore, it will be an interesting study to establish invariance between functionals of Hankel determinants for certain classes of functions.
	Let $\mathcal{A}$ denote the family of all normalized analytic function $f$ of the form \begin{align}\label{e-1.1}
		f(z)=z+\sum_{n=2}^{\infty}a_nz^n,  \;\mbox{for}\; z\in  \mathbb{D}:=\{z:\in\mathbb{C}:|z|<1 \}
	\end{align} 
For $\alpha\in[0,1)$,  $\mathcal{S}^*(\alpha)$ denotes the subclass of $\mathcal{A}$ of functions $f$ such that 
\begin{align}\label{E-1.2}
	\mathcal{S}^*(\alpha)=\bigg\{f\in\mathcal{A}: {\rm Re}\left(\frac{zf^{\prime}(z)}{f(z)}\right)>\alpha,\; z\in\mathbb{D}\bigg\}
\end{align}
called as a starlike function of order $\alpha$. In particular for $\mathcal{S}(0)=: \mathcal{S}^*$, the class of starlike functions. Robertson \cite{Robertson-Ann.Math-1936} introduced the class of starlike functions of order $\alpha$. For $\alpha$=1/2, we get the class $ \mathcal{S}^*(1/2)$ defined as
	\begin{align}\label{e-1.2}
	\mathcal{S}^*(1/2)=\bigg\{f\in\mathcal{A}: {\rm Re}\left(\frac{zf^{\prime}(z)}{f(z)}\right)>1/2,\; z\in\mathbb{D}\bigg\}
	\end{align}
Marx \cite{A.Marx-Math.Ann.-1932} and Strohha\"cker \cite{E.Strohhacker-1933} have proved the following significant result
\begin{align}\label{e-1.3}
		\mathcal{C}\subset \mathcal{S}^*(1/2).
\end{align}
A well-known study in \cite{E. study}, and in  \cite{Duren-1983-NY} reveal that a function $f\in\mathcal{C}$ if, and only if, 
\begin{align*}
{\rm Re}\left(1+\frac{zf^{\prime\prime}(z)}{f^{\prime}(z)}\right)>0,\; z\in\mathbb{D}.
\end{align*}
 Pommerenke \cite{C. Pommerenke-lond-1966} introduced the $ q $-th Hankel determinant for analytic functions $f\in\mathcal{A}$, defined as
\begin{align*}
	H_{q}(n)(f):=\begin{vmatrix}
		a_{n} & a_{n+1} &\cdots& a_{n+q-1}\\ a_{n+1} & a_{n+2} &\cdots& a_{n+q} \\ \vdots & \vdots & \vdots & \vdots \\ a_{n+q-1} & a_{n+q} &\cdots& a_{n+2q-2}
	\end{vmatrix}.
	\end{align*} where  $n\geq 1$ and $q\geq 1$.Finding the rate of growth of the Hankel determinant $|H_q(n)(f)|$ in terms of $ q $ and $ n $ for the whole class $\mathcal{S}\subset \mathcal{A}$ of univalent functions as well as for its sub-classes is a significant problem. Determining the sharp upper bounds of the Hankel determinant is an interesting study and difficult also. Recently, many authors examined the Hankel determinant  $H_{2}(2)(f)=a_2a_4-a_3^2$ and of  $H_{2}(3)(f)=a_3a_5-a_4^2$ order $ 2 $ (see \cite{Kowalczyk-Lecko-Sim-BAMS-2018,N.E. Cho-J.Math-2017,A.Janteng-J.Pure Appl-2006,SK lee-JIP-2013} and references there). The problem of finding sharp bounds of the third Hankel determinant
	\begin{align}\label{e-1.5}
			H_{3}(1)(f):=\begin{vmatrix}
			a_{1} & a_{2} & a_{3}\\ a_{2} & a_{3} & a_{4} \\  a_{3} & a_{4} & a_{5}
		\end{vmatrix}=2a_2a_3a_4-a_3^2-a_4^2+a_3a_5-a_2^2a_5
	\end{align} is a much more difficult task for study. Recently, sharp bounds for $ |H_{3}(1)(f)| $ for the class  $S^*$ of starlike function (see \cite{A.Lecko-Forum-2022}), for the class $\mathcal{C} $ of convex function  (see \cite{Kowalczyk-Lecko-Sim-BAMS-2018}), and for the class of starlike functions of order $ 1/2 $ (see \cite{A.Lecko-CAOT-2019}) have been found. Kowalczyk \emph{et al.} \cite{Kowalczyk-Lecko-Sim-BAMS-2018} have established that $|H_{3}(1)(f)|\leq4/135 $ showing the bound $ 4/135 $ is sharp for $f\in \mathcal{C}$. In view of \eqref{e-1.3} the corresponding bounds of $|H_{3}(1)(f)|$  tell us about the richness of the class. Finding the Sharp bounds for inverse coefficients of the Hankel determinant is a much more difficult task. However, nowadays researchers are engaged in finding the sharp bounds of the Hankel determinant of inverse coefficients of various classes of univalent functions of order $ 2 $ and $ 3 $ (see \cite{Raza-Riaz-Tho-BAMS-2023,Sim-Thomas-Zaprawa-CVEE-2022}). \vspace{1.2mm}

Since the class $\mathcal{S}^*(\alpha)$ subset of $\mathcal{S}$ then by classical result of L\"owner \cite{Lowner-Math Ann-1923} there exists an inverse function $f^{-1}$ in some neighborhood of origin defined on some disc $|w|\leq r_0(f)$ with Taylor series expansion 
\begin{align}\label{e-1.4}
	f^{-1}(w)=w+A_2w^2+A_3w^3+A_4w^4+\cdots.
\end{align} then for $n\geq2,\;|A_n|\leq K_n$, where $K_n=\frac{(2n)!}{n!(n+1)!}$ and $K(w)=w+K_2w_2+K_3w_3+\cdots$ is inverse of the Koebe function, and this inequality is sharp.
We first note that if $f^{-1}$ is the inverse function of $f$, and $f$ and $f^{-1}$ are given by \eqref{e-1.1} and \eqref{e-1.4}, respectively, then comparing coefficients in $f(f^{-1}(w))=w$ gives that
\begin{align}\label{e-2.6}
\begin{cases}
&A_2=-a_2,\\&A_3=2a_2^2-a_3,\\& A_4=5a_2a_3-5a_2^3-a_4,\\& A_5=14a_2^4-21a_3a_2^2+6a_2a_4+3a_3^2-a_5.
\end{cases}
\end{align}
To date, researchers have made progress in discovering a limited number of results regarding the Hankel determinant in relation to inverse functions derived from univalent functions. Recently, Sim \emph{et al.} \cite{Sim-Thomas-Zaprawa-CVEE-2022} have established the sharp bounds of the second Hankel determinant $ |H_{2}(2)(f^{-1})| $ for starlike and convex function of order alpha, whereas Raza \emph{et al.} \cite{Raza-Riaz-Tho-BAMS-2023} found the sharp bound of third Hankel determinant $ |H_{3}(1)(f^{-1})| $ for inverse coefficients of convex functions.\vspace{1.2mm}

 In this paper, we aim to establish invariance between 3rd Hankel determinants $|H_3(1)(f)|$ and $|H_3(1)(f^{-1})|$ for the class $\mathcal{S}^*(1/2)$.\vspace{1.2mm}

To serve our purpose, we find the sharp bound for $|H_3(1)(f^{-1})|$ for $f\in \mathcal{S}^*(1/2)$. It is demonstrated in \cite{Raza-Riaz-Tho-BAMS-2023} that a noninvariance property for $ f\in\mathcal{C} $ between corresponding functionals $ |H_3(1)(f)| $ and $ |H_3(1)(f^{-1})| $ discussed in \cite{Raza-Riaz-Tho-BAMS-2023}, thus settling a conjecture was made in \cite{Kowalczyk-Lecko-Sim-BAMS-2018}. The sharp bound of the third Hankel determinant $ H_3(1)(f^{-1}) $ for the coefficients of the inverse function $ f^{-1} $ of $ f\in\mathcal{C} $ are obtained by  Raza \emph{et al.} \cite{Raza-Riaz-Tho-BAMS-2023}. For the class $ \mathcal{S}^*(1/2) $, Rath \emph{et al.} \cite{Rath-Kumar-Krishna-Lecko-CAOT-2022} provided a revised proof for the result originally established by Lecko \emph{et al.} \cite{A.Lecko-CAOT-2019}. In their study, they determined the sharp bound of the third Hankel determinant $ |H_3(1)(f)| $ and demonstrated that $ |H_3(1)(f)|\leq 1/9 $. \vspace{1.2mm}

Based on the observations discussed above, for further contribution in the topic, it is natural to raise the following questions.
\begin{ques}\label{q-1.1}
	Can we establish the sharp bound of $ |H_3(1)(f^{-1})| $ for the class $\mathcal{S}^*(1/2) $?
\end{ques}
\begin{ques}\label{q-1.2}
	Can we establish the invariance between between the functionals $ |H_3(1)(f)| $ and $ |H_3(1)(f^{-1})| $ for the class $\mathcal{S}^*(1/2) $?
\end{ques}

Before delving into the statement of our result and its accompanying proof, there are some preliminary tasks that need to be addressed. Let $\mathcal{P}$ be the class of analytic functions $p$ defined for $z\in\mathbb{D}$ with Taylor series expansion
	\begin{align}\label{e-1.6}
		p(z)=1+\sum_{n=1}^{\infty}c_nz^n
	\end{align} 
	which have positive real part in $\mathbb{D}$. The following lemma concerning the coefficients of the class $ \mathcal{P} $ will play a key role in proving the main result of the paper.
	\begin{lem}\cite{CMFT-2018,R.J.Libera-American-1982,C. Pommerenke-1975}\label{lem-1.1} Let $p\in\mathcal{P} $ be given by \eqref{e-1.6} with $c_1>0$. Then 	
		\begin{align*}
		&2c_2=c_1^2+\delta(4-c_1^2), \\\;\;&4c_3=c_1^3+2(4-c_1^2)c_1\delta-(4-c_1^2)c_1\delta^2+2(4-c1^2)(1-|\delta|^2)\eta, \\ &8c_4=c_1^4+(4-c_1^2)\delta(c_1^2(\delta^2-3\delta+3)+4\delta)-4(4-c_1^2)(1-|\delta|^2)(c_1(\delta-1)\eta\\&\nonumber\quad+\overline{\delta}\eta^2-(1-|\eta|^2)\rho)
		\end{align*}
		for some $\delta$, $\eta$, and $\rho$ such that $|\delta|\leq1$, $|\eta|\leq1$ and $|\rho|\leq1$.
	\end{lem}
	\section{Main Result}
	\begin{thm}\label{Th-2.1}
		Let $ f\in \mathcal{S}^{*}(1/2) $ be defined by \eqref{e-1.2}. Then
		\begin{align*}
			|H_{3,1}(f^{-1})|\leq\dfrac{1}{9}.
		\end{align*}
		The result is sharp for the function $ f_0\in S^*(1/2) $ given by
		\begin{align}\label{EE-2.1}
			f_0(z)=\frac{z}{(1-z^3)^{1/3}}.
		\end{align}
	\end{thm}
\begin{rem}
The conclusion of Theorem \ref{Th-2.1} reveals the following:
\begin{enumerate}
	\item[(i)] The sharp bound of $|H_3(1)(f^{-1})|$ for the class $\mathcal{S}^*(1/2)$ is $1/9$. Consequently, Question \ref{q-1.1} is answered successfully.
	\item[(ii)] The sharp bound of the functionals $|H_3(1)(f^{-1})|$ and $|H_3(1)(f^{-1})|$ are equal with $1/9$. Hence, the Question \ref{q-1.2} is answered successfully.
\end{enumerate}
\end{rem}
		\begin{proof}[\bf Proof of Theorem \ref{Th-2.1}] Let $f\in\mathcal S^{*}(1/2)$.  Then from (\ref{e-1.2}) there exists a holomorphic  function $p\in\mathcal{P}\;\mbox{of the form} (\ref{e-1.6}) $ such that
			 \begin{align}\label{partha}
				zf^{\prime}(z)=\frac{1}{2}(p(z)+1)f(z), z\in\mathbb{D}
			\end{align} 
			Putting the series of (\ref{e-1.1}) and (\ref{e-1.6}) in (\ref{partha}) and further simplification yeilds the following coefficients
			\begin{align}\label{e-2.2}
				\begin{cases}
					&a_2=\dfrac{1}{2}c_1\vspace{1.2mm}\\&a_3=\dfrac{1}{8}(2c_2+c_1^2)\vspace{1.2mm}\\&a_4=\dfrac{1}{48}(8c_3+6c_1c_2+c_1^3)\vspace{1.2mm}\\&a_5=\dfrac{1}{384}(48c_4+32c_1c_3+12c_2^2+12c_1^2c_2+c_1^4).
				\end{cases}
			\end{align}
			
			Substituting \eqref{e-2.2} into \eqref{e-2.6}, we obtain
			\begin{align}\label{e-2.10}
					&A_2=-\frac{c_1}{2},\\& \label{Eq-2.6} A_3=\frac{3}{8}c_1^2-\frac{1}{4}c_2,\\& \label{Eq-2.7} A_4=-\frac{1}{3}c_1^3+\frac{1}{2}c_1c_2-\frac{1}{6}c_3,\\& \label{Eq-2.8} A_5=\frac{125}{384}c_1^4-\frac{25}{32}c_1^2c_2+\frac{5}{32}c_2^2+\frac{5}{12}c_1c_3-\frac{1}{8}c_4.
			\end{align}
			Then by using \eqref{e-2.10}-\eqref{Eq-2.8} in (\ref{e-1.5}), a simple computation shows that
			\begin{align}\label{e-2.14}
				\nonumber H_3(1)(f^{-1})=\frac{1}{9216}(17c_1^6-102c_1^4c_2+&32c_1^3c_3+180c_1^2c_2^2-144c_1^2c_4+192c_1c_2c_3-216c_2^3\quad\quad\quad&\\+288c_2c_4-256c_3^2)
			\end{align}
		Since the classes $\mathcal{P}$, $\mathcal{S}^*(1/2)$ and the functional $H_3(1)(f^{-1})$ are rotationally invariant, by Carath$ \acute{e} $odory Theorem we may assume, $c:=c_1\in[0,2]$ (see, \cite{C. caratheodory-1907, Duren-1983-NY}) and using the Lemma \ref{lem-1.1} and simplification of \eqref{e-2.14}, we obtain
	
		\begin{align}\label{e-2.15}
			H_3(1)(f^{-1})=\frac{1}{9216}(g_1(c,\delta)+g_2(c,\delta)\eta+g_3(c,\delta)\eta^2)+v(c,\delta,\eta)\rho),
		\end{align} 
	where $ \delta,\eta,\rho\in\overline{\mathbb{D}} $
	and
	\begin{align*}
		\begin{cases}
			g_1(c,\delta):=\quad\delta^2(4-c^2)^2(2c^2-(36-13c^2)\delta)+2c^2\delta^2;\\g_2(c,\delta):=-8c\delta(4-c^2)^2(1+\delta)(1-|\delta|^2);\\g_3(c,\delta):=-8(4-c^2)^2(8+|\delta|^2)(1-|\delta|^2);\\v(c,\delta,\eta):=\quad72\delta(4-c^2)^2(1-|\delta|^2)(1-|\eta|^2).
		\end{cases}
	\end{align*}
			Next, using $|\delta|=x$, $|\eta|=y$ and the fact that $|\rho|\leq1$, we easily obtain 
			\begin{align}\label{e-2.17}
				9216|H_3(1)(f^{-1})|\leq(|g_1(c,\delta)|+|g_2(c,\delta)|+|g_3(c,\delta)|+|v(c,\delta,\eta)|)\leq M(c,x,y),
			\end{align}
		where $ M(c,x,y) $ is defined by
			\begin{align*}
				M(c,x,y):=(h_1(c,x)+h_2(c,x)y+h_3(c,x)y^2+h_4(c,x)(1-y^2))
			\end{align*}
		and
		\begin{align*}
			\begin{cases}
			h_1(c,x):=&x^2(4-c^2)^2(2c^2+(36-13c^2)x)+2c^2x^2;\\	h_2(c,x):=&8cx(4-c^2)^2(1+x)(1-x^2);\\	h_3(c,x):=&8(4-c^2)^2(8+x^2)(1-x^2);\\	h_4(c,x):=&72x(4-c^2)^2(1-x^2).
			\end{cases}
		\end{align*}
		The goal is to find the maximum values of $M(c,x,y)$ on the closed cuboid $\Omega=[0,2]\times[0,1]\times[0,1]$. This involves determining the maximum values in the interior of $\Omega$, the interior of the six faces, as well as at the vertices and edges (twelve in total).\vspace{1.2mm}
		
		\noindent{\bf Case 1.}	\textit{In the interior of $\Omega$:}\\ A simple computation shows that 
		\begin{align*}
				\frac{\partial M(c,x,y)}{\partial y} =0\;\mbox{implies that}\; 8(4-c^2)^2(1-x^2)[cx(1+x)+2(1-x)(8-x)y]=0
			\end{align*}
		only for 
		\begin{align*}
			y=-\frac{cx(1+x)}{2(8-x)(1-x)}=y_0(c,x)
		\end{align*}
			since  $(4-c^2)^2\neq 0 $ and $(1-x^2)\neq0 $ for all $c\in(0,2)$ and $x\in (0,1) $. Moreover, we see that $y_0(c,x)<0$ for all $c\in(0,2)$ and $x\in (0,1) $. Hence, we deduce that $M$ has no critical points in interior of $\Omega$.\vspace{1.2mm}
			
	\noindent{\bf Case 2.} \textit{On the edges of $ \Omega$.}\\
		\noindent{\bf Case 2.1.} On $ x=1, y=0  $ or $ x=1,y=1 $, we see that
		\begin{align*}
			M(c,1,0)=M(c,1,1)=(4-c^2)^2(36-13c^2)\leq 576,\; c\in(0,2).
		\end{align*}
	\noindent{\bf Case 2.1.} On	$ x=0, y=1,$ we have
	\begin{align*}
	M(c,0,1)=64(4-c^2)^2\leq1024,\;c\in(0,2),
		\end{align*}
		which is equivalent to 
		\begin{align*}
			64c^4-512c^2\leq 0\; \mbox{for}\;c\in(0,2)
		\end{align*}
	which is true.\vspace{1.2mm}
	
	\noindent{\bf Case 2.3.} If $c=0,y=0$, then an easy computation shows that
	 \begin{align*}
		M(0,x,0)=-576x^3+1152x\leq M(0,x_1,0)\leq256\sqrt{6}\approx627.069,\; x\in(0,1)
	\end{align*}
			where\begin{align*}
				x_1=\sqrt{\frac{2}{3}}\approx0.816497.			\end{align*}
		\noindent{\bf Case 2.4.}  If $ c=0,y=1 $, then we see that \begin{align*}
			M(0,x,1)=1024-896x^2+576x^3-128x^4\leq1024,\;x\in(0,1)
			\end{align*}
		which is equivalent to 
		\begin{align*}
		-128x^2+576x-896\leq 0,\;x\in(0,1)
			\end{align*}
		which is true.\vspace{1.2mm}
		
		\noindent{\bf Case 2.5.} If $ c=0,x=0 $, then 
		we see that
		\begin{align*}
			\\M(0,0,y)=1024y^2\leq1024,y\in(0,1).
			\end{align*}
		\noindent{\bf Case 2.6.} If $ c=0, x=1, $ then
		\begin{align*}
			M(0, 1, y)=576,\;y\in(0,1).
			\end{align*}
			\noindent{\bf Case 2.7.} For the other edges,  $ c=2,x=0  $or $ c=2,x=1 $ or $ c=2,y=0 $ or $ c=2,y=1 $ or $ x=0,y=0 $, we have
			\begin{align*}
				M(2,0,y)=M(2,1,y)=M(2,x,0)=M(2,x,1)=M(c,0,0)=0
			\end{align*}
		for all $ c \in(0,2) $, $ x\in(0,1) $ and $ y\in(0,1) $.\vspace{1.2mm}
		
			\noindent{\bf Case 3.} \textit{On the faces of $\Omega$.}\vspace{1.2mm}
			
			 \noindent{\bf Case 3.1.} On $ c=2 $, we have \begin{align*}
				M(2,x,y)=0, x\in(0,1), y\in(0,1).
			\end{align*}
			\noindent{\bf Case 3.2.} On $ c=0 $, we see that \begin{align*}
				M(0,x,y)&=1152-576x^3+(1024-1152x-896x^2+1152x^3-128x^4)y^2\\&=576(2-x^2)+128(8-x)(1-x)^2(1+x)y^2\\&\leq576(2-x^2)+128(8-x)(1-x)^2(1+x)\\&=1024-896x^2+576x^3-128x^4\;\leq1024,(x,y)\in(0,1)\times(0,1)
			\end{align*}
			which is equivalent to 
			\begin{align*}
				-896x^2+576x^3-128x^4\leq 0 \;\mbox{for all}\;x\in(0,1)
			\end{align*}
			which is true.\vspace{1.2mm}
			
			\noindent{\bf Case 3.3.} On $ x=0, $ we have
			\begin{align*}
				M(c,0,y)=64(4-c^2)^2y^2\leq1024,c\in(0,2),y\in(0,1).
			\end{align*}
		\noindent{\bf Case 3.4.}
		On $ x=1, $ we have
		\begin{align*}
			M(c,1,y)=(4-c^2)^2(36-c^2)\leq576,\;c\in(0,2),y\in(0,1).
			\end{align*}
		\noindent{\bf Case 3.5.} On $ y=0 $, we have
		\begin{align*}
			M(c,x,0)=(4-c^2)^2\left(72x(1-x^2)+x^2(2c^2+(36-13c^2)x+2c^2x^2)\right).
			\end{align*}
		Then a simple computation shows that 
		\begin{align*}
		\quad\quad\quad\quad\frac{\partial }{\partial c}M(c,x,0)=2c(4-c^2)((8-6c^2)x^4+(39c^2+20)x^3+(8-6c^2)x^2-144x)=0
		\end{align*}
		and
		\begin{align*}
         \frac{\partial}{\partial x} M(c,x,0)=(4-c^2)^2(72+4c^2x-(108+39c^2)x^2+8c^2x^3)=0.
			\end{align*}The above system of equation has no critical point on $(0,2)\times(0,1)$, and\begin{align*}
			\max\{M(c,x,0)\}\leq256\sqrt{6}\;\; \mbox{for}\;  {c=0,x=\sqrt{{2}/{3}}}.
			\end{align*} \vspace{1.1mm} 
		
		\noindent{\bf Case 3.6.} On $ y=1 $, we have
		\begin{align*}
			M(c,x,1)=(4-c^2)^2[x^2(2c^2x^2+(36-13c^2)x+2c^2)+8cx(1+x)(1-x^2)+8(8+x^2)(1-x^2)].
			\end{align*}
		Furthermore, we see that
		\begin{align*}
				\quad\quad\quad\frac{\partial}{\partial c}M(c,x,1)=&(4-c^2)[32x(1-x)(1+x)^2-40c^2x(1-x)(1+x)^2\\&-6c^3x^2(2-13x+2x^2)+8c(-32+30x^2-31x^3+6x^4)]=0
			\end{align*}
		and
		\begin{align*}
			\frac{\partial}{\partial x}M(c,x,1)&=(4-c^2)^2[c^2x(4-39x+8x^2)-4x(28-27x+8x^2)\\&\quad-8c(-1-2x+3x^2+4x^3)]=0.
			\end{align*}
		However, it can be shown that the above system has no solutions on $(0,2)\times(0,1)$. Thus, it follows that $ M(c, x, 1) $ has no critical point on $(0,2)\times(0,1)$ and hence,
		\begin{align*}
			\max M(c,x,1)=1024
			\end{align*} 
		attained at $ c=0,x=0 $.\vspace{1.2mm}
		
		\noindent{\bf Case 4.} \textit{On the vertices of $\Omega$}. A simple computation gives the following estimates
		\begin{align*}
			&M(0,0,0)=M(2,0,0)=M(2,1,0)=M(2,1,1)=M(2,0,1)=0,\\	&M(0,1,0)=M(0,1,1)=576,	M(0,0,1)=1024.
			\end{align*}
		Summerizing Case 1 to Case 4, we easily obtain that \begin{align*}
			\max \{M(c,x,y):(c,x,y)\in \Omega\}\leq 1024
			\end{align*}
		Consequently, in view of inequality \eqref{e-2.17}, we obtain that
		 \begin{align*}
				|H_3(1)(f^{-1})|\leq\dfrac{1}{9}.			\end{align*}
		Thus the bound of the theorem is established.\vspace{1.2mm}
		
		To complete the proof, it is sufficient to show that this bound is sharp. In order to show that we consider the function given by \eqref{EE-2.1}, which is equivalent to choosing  $a_2=a_3=a_5=0$ and $a_4=1/3$. By \eqref{e-2.6}, a simple computation shows that $A_2=A_3=A_5=0$ and $A_4=1/3$. Thus it follows from \eqref{e-1.5} that  $|H_3(1)(f^{-1})|=1/9$ establishes that the bound is sharp. This completes the proof.
		\end{proof}
	

		\noindent\textbf{Compliance of Ethical Standards}\\
		
		\noindent\textbf{Conflict of interest} The authors declare that there is no conflict of interest regarding the publication of this paper.\vspace{1.2mm}
		
		\noindent\textbf{Data availability statement}  Data sharing not applicable to this article as no datasets were generated or analyzed during the current study.

\end{document}